\documentclass{amsart}
\usepackage{graphics}
\usepackage{color}
\usepackage{amsmath}
\usepackage{amsthm}
\usepackage{amsfonts}
\newtheorem{theorem}{Theorem}

\newtheorem{corollary}[theorem]{Corollary}
\newcommand{\baseRing}[1]{\ensuremath{\mathbb{#1}}}
\newcommand{\N}{\baseRing{N}}

\newcommand{\R}{\baseRing{R}}

\newcommand{\C}{\baseRing{C}}

\begin{document}
\title[Polynomial Graphs With Applications To Game Theory]{Polynomial Graphs With Applications\\ To Graphical Games, Extensive-Form
Games,\\ and Games With Emergent Node Tree Structures}
\author{Ruchira S. Datta}
\email{datta@math.berkeley.edu}
\urladdr{http://math.berkeley.edu/\~datta}
\date{\today}
\begin{abstract}
We prove a theorem computing the number of solutions to a
system of equations which is generic subject to the sparsity conditions
embodied in a graph.  We apply this theorem to games obeying
graphical models and to extensive-form games.  We define 
\emph{emergent-node tree structures} as additional structures which 
normal form games may have.  We apply our theorem to games having
such structures.  We briefly discuss how emergent node tree structures
relate to cooperative games.

\end{abstract}
\maketitle
The set of Nash equilibria for a game with generic payoff functions is
finite \cite{harsanyi}.  This implies that the set of totally mixed Nash
equilibria for a game with generic payoff functions is also finite.  These
are the real solutions to a system of polynomial equations and
inequalities.  The complex solutions to the system of equations are called
\emph{quasiequilibria}.  Thus, the set of totally mixed Nash equilibria is
a subset of the set of quasiequilibria.  In fact, the set of
quasiequilibria is also finite in the most generic case, and its
cardinality can be computed as a function of the numbers of pure strategies
of the players.  Thus, this is an upper bound on the number of totally
mixed Nash equilibria.  Even in a nongeneric case, as long as the set of
quasiequilibria is finite, its cardinality will be bounded above by the
number in the generic case.

For the main theorem of this article, Theorem \ref{nashbernstein}, we
hypothesize a set of technical conditions that a system of polynomial
equations may satisfy, which are encoded in an associated graph, the
\emph{polynomial graph}, and we
prove a formula describing the number of solutions in this case.  We then
show how to associate such a graph to three special classes of games.  The
first two are graphical games and extensive-form games.  The last is games
with \emph{emergent node tree structure}, a new model for games in which
the players can be hierarchically decomposed into groups.  Usually such
hierarchical decomposition is modelled by \emph{cooperative games}, and we
briefly discuss how our model is related to, yet differs from, the
cooperative framework.

\section{Generic Number of Roots of a Sparse Polynomial System}


The following theorem tells us the number of 0-dimensional complex roots (none
of whose components are zero) of a system of polynomial equations which obeys
certain sparsity conditions and is otherwise generic.  Our formulation of this
theorem is motivated by the applications to game theory which follow, although
such polynomial systems may arise in other contexts.

\begin{theorem}
\label{nashbernstein}
Suppose that $0<d\in\N$ and that we are given a partition
$\{1,\ldots,d\}=\coprod_{i=1}^N T_i$ of
$\{1,\ldots,d\}$.  Write $d_i=|T_i|$.  Suppose further
that we are given a directed graph $G$, the \emph{polynomial graph}, on $d$
vertices, denoted $v_1,\ldots,v_d$, without self-loops and with the
property that for any $v_j$ and $T_i$, if there is some $k\in T_i$ such
that there is an edge from $v_j$ to $v_k$ in $G$, then for every $k\in T_i$
there is an edge from $v_j$ to $v_k$ in $G$.  Let
\begin{eqnarray*}
\label{graphicalmodelgameeqns}
f_1(\sigma_1,\ldots,\sigma_d)&=0,\\
f_2(\sigma_1,\ldots,\sigma_d)&=0,\\
\vdots\\
f_d(\sigma_1,\ldots,\sigma_d)&=0\\
\end{eqnarray*}
be a system (\ref{graphicalmodelgameeqns}) of $d$ polynomial equations in $d$ variables
$\sigma_1,\ldots,\sigma_d$ with the following properties:
\begin{enumerate}
\item All monomials occuring in the $f_i$'s are squarefree.
\item If $\sigma_j,\sigma_k\in T_i$ with $j\neq k$ then $\sigma_j$ and $\sigma_k$ do not
both occur in any monomial of any of the $f_i$'s.
\item If there is no edge from $v_j$ to $v_k$ in $G$ then the variable
$\sigma_k$ does not occur in $f_j$.
\end{enumerate}
Thus, the equations are multilinear, and they are linear over the variables
from each $T_i$.
Construct a $d\times d$ matrix $M$ as follows:  If variable $\sigma_k$ occurs in
the polynomial $f_j$, with $T_i$ the subset containing $v_k$, then
$$M_{jk}={1\over(d_i!)^{1/d_i}},$$ otherwise $M_{jk}=0$.  If the 
system (\ref{graphicalmodelgameeqns}) is $0$-dimensional, then the number
of its solutions in $(\C^*)^d$ (i.e. such that $\sigma_k\neq 0$ for all
$k$) is bounded above by the permanent of $M$, and is equal to the
permanent of $M$ for generic coefficients.
\end{theorem}
\begin{proof}
Without loss of generality, assume
$$T_i=\left\{1+\sum_{l=1}^{i-1}d_l,2+\sum_{l=1}^{i-1}d_l,\ldots,
d_i+\sum_{l=1}^{i-1}d_l\right\},$$
that is, that the $T_i$'s are contiguous.  

Let $a_{ij}=1$ if there is an edge in $G$ from $v_j$ to $v_k$ for $k\in
T_i$, and $a_{ij}=0$ otherwise.  Then the Newton polytope $P_j$ of
$f_j$ is the Cartesian product $P_{1j}\times P_{2j}\times\cdots\times
P_{Nj}$, where $P_{ij}$ is the convex hull of the scaled coordinate vectors
$\left\{a_{ij}e_k\mid k\in T_i\right\}$ and the origin.  For $i$ with
$a_{ij}=1$, $P_{ij}$ is the $d_i$-dimensional unit simplex, and for $i$ with
$a_{ij}=0$, $P_{ij}$ degenerates to the $d_i$-dimensional origin (which is a
$0$-dimensional simplex).  By the Bernstein-Kouchnirenko Theorem
\cite{Bernstein} \cite{Kouchnirenko}, it
suffices to show that the mixed volume of the polytopes $P_1,\ldots,P_d$ is
given by the permanent of $M$.

Let $Q_j=\lambda_1 P_1+\cdots+\lambda_j P_j$, where $+$ denotes Minkowski
addition and the scale factors $\lambda_1,\ldots,\lambda_j$ are parameters.
We show by induction on $j$ that $Q_j=Q_{1j}\times Q_{2j}\times
\cdots\times Q_{Nj}$, where $Q_{ij}$ is the convex hull of
$$\left\{(a_{i1}\lambda_1+a_{i2}\lambda_2+\cdots+a_{ij}\lambda_j)e_k\mid
k\in T_i\right\}$$ and the origin.  (If
$a_{i1}\lambda_1+a_{i2}\lambda_2+\cdots+a_{ij}\lambda_j=0$ then $Q_{ij}$
degenerates to the origin.)  The base case follows from our
characterization of $P_j$ above.  Now consider the Minkowski sum of
$Q_j=Q_{1j}\times\cdots\times Q_{Nj}$ and $\lambda_{j+1}
P_{j+1}=(\lambda_{j+1} P_{1(j+1)})\times\cdots\times(\lambda_{j+1}
P_{N(j+1)})$.  It follows from the definition of Minkowski sum that this
is $(Q_{1j}+\lambda_{j+1}
P_{1(j+1)})\times\cdots\times(Q_{Nj}+\lambda_{j+1} P_{N(j+1)})$, and (using
the induction hypothesis) that each factor $Q_{ij}+\lambda_{j+1}P_{i(j+1)}$
is equal to the convex hull of
$$\left\{(a_{i1}\lambda_1+a_{i2}\lambda_2+\cdots+a_{ij}\lambda_j+a_{i(j+1)}e_{j+1})e_k\mid
k\in T_i\right\}$$ and the origin.  

The $d_i$-dimensional volume of the $d_i$-dimensional unit simplex scaled
by $\lambda$ in each dimension is
$${\lambda^{d_i}\over (d_i)!}.$$
We are interested in the $d$-dimensional volume of $Q_d$.  If
$a_{i1}=a_{i2}=\cdots=a_{id}=0$ for some $i$, then this volume vanishes,
and hence the mixed volume also vanishes.  In this case the $k$th column of the
matrix $M$ will be all zeroes for any $k\in T_i$, so the permanent of $M$
also vanishes, and the theorem holds.  So assume that for each $i$, there
is some $j$ with $a_{ij}=1$.  Then the volume of $Q_d$ is
$$\prod_{i=1}^N
{\left(a_{i1}\lambda_1+\cdots+a_{id}\lambda_d\right)^{d_i}
\over d_i!}.$$
Let $(g_{jk})$ be the adjacency matrix of $G$, that is, $g_{jk}=1$ if there
is an edge in $G$ from $v_j$ to $v_k$ and $g_{jk}=0$ otherwise.  Then
$a_{ij}=g_{jk}$ for all $k\in T_i$.  So the volume of $Q_d$ is
$${\prod_{k=1}^d \left(g_{1k}\lambda_1+\cdots+g_{dk}\lambda_d\right)
\over\prod_{i=1}^N d_i!}.$$
The mixed volume of $P_1,\ldots,P_d$ is the coefficient of
$\lambda_1\lambda_2\cdots\lambda_d$ in the above expression, which is the
permanent of $(g_{jk})$ divided by $\prod_{i=1}^N d_i!$.

It remains to show that the permanent of $M$ is the permanent of $(g_{jk})$
divided by $\prod_{i=1}^N d_i!$.  Note that $M_{jk}\neq 0$ exactly when
$g_{jk}\neq 0$.  We induct on $N$.  For the base case,
$d_1=d$, and each nonzero entry of $M$ is $(1/d!)^{1/d}$.  A term in the
permanent of $M$ is the product of $d$ entries from $M$, so if it is
nonzero it is $1/d!$.  Thus the permanent of $M$ is $1/d!$ times the
permanent of $(g_{jk})$, as required.  Now partition the matrix $M$ and the
matrix $(g_{jk})$ into two vertical bands corresponding to the subsets
$\cup_{i=1}^{N-1} T_i$ and $T_N$.  The permanent can be computed as
the sum of a term for each choice of $d_N$ rows $1\leq j_1<\cdots<j_{d_N}\leq
d$: compute the $(d-d_N)\times(d-d_N)$ subpermanent of the left band
obtained by crossing out those rows, compute the $d_N\times d_N$
subpermanent of the right band corresponding to those rows, and multiply
them together.  By the inductive hypothesis, the left subpermanent of $M$
is the left subpermanent of $(g_{jk})$ divided by $\prod_{i=1}^{N-1} d_i!$.  
For the right subpermanent, every row is either all nonzero or all zero.
If any row is all zero, both right subpermanents vanish.  If every entry is
nonzero, then all the entries are the same: $g_{jk}=1$ and
$M_{jk}=(1/d_N!)^{1/d_N}$.  The right subpermanent of $M$ is
$d_N!\left(\left(1/d_N!\right)^{1/d_N}\right)^{d_N}=1$, and the right
subpermanent of $(g_{jk})$ is $d_N!$.  So the whole term for $M$ is the
whole term for $(g_{jk})$ divided by $\prod_{i=1}^N d_i!$.
\end{proof}

We note that if the coefficients are generic subject to the conditions
given in Theorem \ref{nashbernstein}, {\em all} the solutions to the system
will lie in the torus $(\C^*)^d$.  In what follows we will refer to ``the
number of solutions in the torus $(\C^*)^d$'' as ``the number of
solutions'' by abuse of language.

\begin{corollary}
Convert the directed graph $G$ of Theorem \ref{nashbernstein} into a bipartite
graph on $2d$ vertices, with the source of every edge on the left side and
the target of every edge on the right side.
If the system in Theorem \ref{nashbernstein} is $0$-dimensional with
generic coefficients, then it
has a solution if and only if this bipartite graph has a perfect matching.
\end{corollary}
\begin{proof}
From the proof of Theorem \ref{nashbernstein}, we see that the number of
solutions is nonzero if and only if the permanent of the adjacency matrix
is nonzero.  It is a well-known fact that this is equivalent to the
existence of a perfect matching: any permutation $\pi$ which contributes a
nonvanishing term $\prod_{j=1}^d g_{j\pi(j)}$ to the permanent corresponds
to a perfect matching, where vertex $j$ on the left is matched to vertex
$\pi(j)$ on the right.
\end{proof}

In fact, we could have used the bipartite graph in Theorem
\ref{nashbernstein}. However, we defined the polynomial graph to be the
directed graph to remain consistent with the usual definition of graphical
models of games.  

\begin{corollary}
\label{cyclicgraph}
If the system in Theorem \ref{nashbernstein} is $0$-dimensional and has a
solution, then every node in the graph $G$ lies on a directed cycle.
\end{corollary}
\begin{proof}
As in the proof of the previous corollary, a permutation $\pi$ must exist
such that $j$ has an edge to $\pi(j)$ for every $j=1,\ldots,d$.
This permutation can be expressed as a product of disjoint cycles.  Each
node lies in one of these cycles, and a cycle of the permutation
corresponds to a directed cycle in the graph.
\end{proof}

We should note carefully that the 
Bernstein-Kouchnirenko theorem gives the number of solutions to a {\em
0-dimensional} polynomial system.  So when the number given by that
theorem---in particular, the permanent of the matrix in Theorem
\ref{nashbernstein}---vanishes, {\em either} the polynomial system has no
solution, {\em or} its solution set has positive dimension.

Note that the conditions on $G$
imply that the matrix $M$ has a $d_i\times d_i$ block of zeroes along its
diagonal for $i=1,\ldots,N$.  This is because $G$ has no self-loops, and if
it had an edge from an element $v_j$ of $T_i$ to any other element $v_k$ of
$T_i$, then there would have to be an edge from $v_j$ to every element of
$T_i$ including itself.

\section{Finite Games}
\noindent
In the remainder of this article, we apply Theorem \ref{nashbernstein} to game
theory in a few different contexts.  We now introduce the notation 
we will need from game theory.  The concepts we describe in this section can be found in a standard game
theory text such as \cite{OsborneRubinstein}. However, in some cases we use
simplified notation for the restricted situations we will consider.

Game theory is the study of strategic interaction.
Such interaction takes place between multiple agents in a single setting,
or \emph{environment}.  An \emph{agent} is an entity which can receive
\emph{information} about the state of the environment (including itself and
other agents), take \emph{actions} which may alter that state, and express
\emph{preferences} among the various possible states.  These preferences
are encoded for each agent by a \emph{utility function}, a mapping from the
set of all states to $\R$.  Its value for a particular state is the
\emph{utility} of that state for the agent.  The agent prefers one state to
another if its utility is greater, and is indifferent between them if
their utilities are equal.\footnote{%
Instead of specifying the utility of each state for each agent, one might
specify the {\em change\/} in utility, or \emph{marginal utility}, which
accrues to each agent upon each {\em transition\/} between states.  Clearly
any utility function induces a marginal utility function, but unless one
imposes additional conditions a marginal utility function may not
induce a utility function.  Such a marginal utility function, which one
might call \emph{intransitive}, could still be a useful model of reality.
For example, one wouldn't necessarily feel the same about being laid off
and then immediately rehired as if one had simply continued in the same
position.  However, we will not consider such intransitive marginal utility
functions any further.
}
Changes in the state of the environment may
also occur spontaneously (i.e., not due to the actions of any of the
agents).  A \emph{strategy} is a (possibly stochastic) rule for
an agent to choose an action at every point when the agent may act, given
the available information.  A \emph{rational} agent is one whose strategy
maximizes its expected utility under the circumstances.

We will restrict attention to games which take place in a finite number of
time steps between a finite number of agents, each of which has a finite
number of possible actions.  The agents are called \emph{players}, and
whenever they take an action they are said to \emph{move}.  A spontaneous
change in the state of the environment is called a \emph{move by nature}.
The game is over when no player (including nature) has any possible
actions.  The state of the environment at such a terminal stage is called
an \emph{outcome}.  Generally preferences are specified only over outcomes,
not at intermediate stages of the game.

The first type of game we will consider is the \emph{normal-form game}.  In
a normal-form game, there is only one time step, at which all the players
move simultaneously.  We denote the set of players by $I=\{1,\ldots,N\}$.
The actions a player can take are called \emph{pure strategies}.  We associate
to the players finite disjoint sets of
pure strategies $S_1,\ldots,S_N$.  For each $i$ let $d_i=|S_i|-1$.  We
write the set $S_i$ as $\{s_{i0},\ldots,s_{id_i}\}$.  We write
$S=\prod_{i\in I}S_i$.  Game play consists of the collective choice of an
element of $S$ by the players: each player $i$ moves by choosing an element
of $S_i$.  We identify $S$ as the set of possible outcomes.  We denote by
$u_i(s)$ the utility for player $i$ of the outcome $s\in S$.  Thus, the
game is completely specified by the number $N$ of players, the sets $S_i$
of pure strategies, and the utility functions (or \emph{payoff functions})
$u_i\colon S\mapsto\R$.

A player may move stochastically rather than deterministically.  In that
case the player is said to execute a \emph{mixed strategy}.  The mixed
strategy specifies the probability with which the player chooses each
possible action.  The set $\Sigma_i$ of mixed strategies of player
$i$ is the set of all functions $\sigma_i\colon S_i\mapsto [0,1]$ with
$\sum_{s_{ij}\in S_i}\sigma_i(s_{ij})=1$.  That is, it is the
$d_i$-dimensional probability simplex.  We write $\Sigma=\prod_{i\in
I}\Sigma_i$.  An element $\sigma$ of $\Sigma$, which specifies the
strategies executed by all the players, is called a \emph{strategy
profile}.  If the players execute the strategy profile $\sigma$, then the
probability of outcome $s$ is $\sigma(s)=\prod_{i=1}^N \sigma_i(s_i)$.  The
\emph{expected} utility for player $i$ of the strategy profile $\sigma$ is
given by multilinearity as $u_i(\sigma)=\sum_{s\in S}u_i(s)\sigma(s)$.

When considering how agent $i$ should behave, it will be convenient to
separate out $i$'s own strategy, over which $i$ has control, from the
strategies of all the other players.  We write $\Sigma_{-i}=\prod_{j\in
I-\{i\}}\Sigma_j$, and we write $\sigma_{-i}$ for the image of
$\sigma\in\Sigma$ under the projection $\pi_{-i}$ from $\Sigma$ onto
$\Sigma_{-i}$.  By abuse of notation, we write $u_i(\tau_i,\sigma_{-i})$ for
the $i$th player's expected payoff from the strategy $\sigma$ whose $i$th
component is $\tau_i$ and whose other components are defined by
$\pi_{-i}(\sigma)=\sigma_{-i}$.

We assume \emph{perfect information}: each player knows the complete
specification of the game, knows that every player knows, knows that
every player knows that every player knows, ad infinitum.  That is, the
specification of the game is \emph{common knowledge}.  Under these
circumstances, what is rational behavior?  In his landmark paper
\cite{nash:definition}, John Nash answered this question in terms of what is
now called \emph{best response}.  A best response of player $i$ to the
strategy profile $\sigma$ is a mixed strategy $\sigma_i^*$ such that
$u_i(\sigma_i^*,\sigma_{-i})\geq u_i(\sigma'_i,\sigma_{-i})$ for any other
mixed strategy $\sigma'_i$ of player $i$.  That is, given that all the
other players execute the strategy profile $\sigma_{-i}$, the mixed
strategy $\sigma_i^*$ maximizes player $i$'s expected utility.  A
\emph{Nash equilibrium} is a strategy profile which is a best response to
itself for all the players.  That is, it is a strategy profile $\sigma^*$
such that for each player $i$, we have $u_i(\sigma^*)\geq
u_i(\sigma'_i,\sigma^*_{-i})$ for every other mixed strategy $\sigma'_i$ of
player $i$.  Nash proved that such an equilibrium always exists.

How can we compute the Nash equilibria of a given game?
We need to search the set $\Sigma$ of strategy profiles, which is a
polytope: the product of probability simplices.  We can decompose the
problem by {\em stratifying\/} this polytope: first we look for Nash
equilibria in its interior, then in the interiors of its facets, then in
the interiors of the facets of those facets, and so forth, until finally we
look for Nash equilibria at the vertices of the polytope (that is, pure
strategy Nash equilibria).  A strategy profile $\sigma$ lies in the
interior of this polytope if $\sigma_i(s_{ij})>0$ for every $s_{ij}\in
S_i$, for every $i$.  Such a strategy profile is called \emph{totally
mixed}.  Note that a totally mixed Nash equilibrium {\em need not} exist.

So we concentrate our attention on the totally mixed Nash equilibria.  We
observe that for a totally mixed strategy profile $\sigma$ to be a Nash
equilibrium, it is necessary and sufficient that for each player $i$ we
have 
$u_i(s_{ij},\sigma_{-i})=u_i(s_{ik},\sigma_{-i})$ for any pure strategies
$s_{ij},s_{ik}\in\Sigma_i$.  These equations are called the
\emph{indifference equations} for player $i$.  The sufficiency is clear.
For the necessity,
suppose to the contrary that
$u_i(s_{ij},\sigma_{-i})>u_i(s_{ik},\sigma_{-i})$.  Define $\sigma'_i$ by
$$\sigma'_i(s_{il})=
\begin{cases}
\sigma_i(s_{ij})+\sigma_i(s_{ik}),&l=j\\
0,&l=k\\
\sigma_i(s_{il}),&\hbox{otherwise}\\
\end{cases}.$$
Then since $\sigma_i(s_{ik})>0$, we have
$$u_i(\sigma_i',\sigma_{-i})=u_i(\sigma)+\sigma_i(s_{ik})\left(u_i(s_{ij},\sigma_{-i})-u_i(s_{ik},\sigma_{-i})\right)>u_i(\sigma),$$
a contradiction.  

So we have a system of $\sum_{i=1}^N d_i$ equations,
$u_i(s_{ij},\sigma_{-i})=u_i(s_{i0},\sigma_{-i})$ for $j=1,\ldots,d_i$, for
$i=1,\ldots,N$, in $\sum_{i=1}^N d_i$ unknowns $\sigma_i(s_{ij})$ for
$j=1,\ldots,d_i$, for $i=1,\ldots,N$. (Here we have \emph{dehomogenized},
that is, we have eliminated $\sigma_i(s_{i0})$ by substituting
$1-\sum_{j=1}^{d_i}\sigma_i(s_{ij})$).  What we are equating are the
expressions $u_i(s_{ij},\sigma_{-i})=\sum_{s_{-i}\in
S_{-i}}u_i(s_{ij},s_{-i})\sigma_1(s_1)\cdots\sigma_{i-1}(s_{i-1})\sigma_{i+1}(s_{i+1})\cdots\sigma_N(s_N)$,
which are multilinear polynomials whose coefficients are the real numbers
$u_i(s)$.  The (possibly complex) roots of this system are called
\emph{quasiequilibria}, and those roots which are totally mixed strategy
profiles (that is, which are real with $\sigma_i(s_{ij})>0$ and
$\sum_{j=1}^{d_i}\sigma_i(s_{ij})<1$) are the totally mixed Nash
equilibria.  

Now we see how Theorem \ref{nashbernstein} applies to normal-form games.  In
this case, each $T_i$ corresponds to the
set of strategies of player $i$.  The blocks of zeroes along the diagonal
imply that a player's expected payoffs from their own pure strategies do
not depend on the probabilities they have assigned to their own pure
strategies, so these polynomial systems do indeed correspond to the
equations for totally mixed Nash equilibria of games. 
\begin{corollary}
\label{beitragetheorem}
Consider a normal form game between players $I=\{1,\ldots,N\}$ with pure
strategy sets $S_i$ for each $i$ and generic utility functions $u_i\colon
\prod_{i\in I}S_i\to\R$.  Construct a graph $G$ with nodes $\coprod_{i\in
I}(S_i-\{s_{i0}\})$
such that there is an edge from $s_{ik}$ to $s_{jl}$ in $G$ if and only if
$i\neq j$.  Let the variable corresponding to $s_{ik}$ be
$\sigma_i(s_{ik})$ and the equation corresponding to $s_{ik}$ be the
indifference equation $u_i(s_{ik},\sigma_{-i})=u_i(s_{i0},\sigma_{-i})$.
Then this system of equations obeys the conditions of Theorem
\ref{nashbernstein}, so the number of solutions in the generic case is
given by that theorem.
\end{corollary}
This special case was proved as Theorem 2 in \cite{mclennan}, so our Theorem
\ref{nashbernstein} is a generalization of that theorem.

\section{Graphical Games}

Kearns, Littman, and Singh
\cite{kearnslittmansingh} defined the
concept of {\em graphical games}, or games obeying \emph{graphical models}.
(That paper considers undirected graphs, but the extension to directed
graphs which we will use is straightforward.)  A game between players
$1,\ldots,N$ obeys a directed graphical model, if the payoffs to player
$i_1$
only depend on the actions of those players $i_2\neq i_1$ for which there is an
edge from $i_1$ to $i_2$ in the graphical model.

Our theorem applies
in particular to graphical games. As in Corollary
\ref{beitragetheorem}, we take the pure strategy sets $S_i$ to be the sets
$T_i$ of Theorem \ref{nashbernstein}.  Given a
polynomial graph
$G$ as in Theorem \ref{nashbernstein}, we draw an edge
from $i_1$ to $i_2$ in the graphical model if there is any $j\in T_{i_1}$
with edges to the vertices
in $T_{i_2}$ in the polynomial graph $G$.  The polynomial graph $G$ may not
represent the most generic 
case of the graphical model, however.  If we are given a graphical model,
then to construct its polynomial graph $G$, for any edge from $i_1$ to
$i_2$, we draw edges in $G$ from {\em every} vertex $j\in T_{i_1}$ to every
vertex in $T_{i_2}$.  

\begin{corollary}
\label{graphicalmodeltheorem}
Suppose a normal form game between players $I=1,\ldots,N$ with pure
strategy sets $S_i$ for each $i$ and utility functions $u_i\colon
\prod_{i\in I}S_i\to\R$ obeys a directed graphical model $\gamma$ with
nodes $1,\ldots,N$.  Construct a graph $G$ with nodes $\coprod_{i\in I}S_i$
such that there is an edge from $s_{ik}$ to $s_{jl}$ in $G$ if and only if
there is an edge from $i$ to $j$ in $\gamma$.  Then the system of equations
defining the quasiequilibria of $G$ satisfies the hypotheses
of Theorem \ref{nashbernstein}, so the number of such quasiequilibria in
the generic case is
given by the permanental formula.
\end{corollary}

For example, consider a game with 4 players, each with 3 pure strategies.
Generically, such a game has
\begin{equation*}
\hbox{per}
\begin{pmatrix}
0&0&\frac{1}{\sqrt2}&\frac{1}{\sqrt2}&
\frac{1}{\sqrt2}&\frac{1}{\sqrt2}&\frac{1}{\sqrt2}&\frac{1}{\sqrt2}\\
0&0&\frac{1}{\sqrt2}&\frac{1}{\sqrt2}&
\frac{1}{\sqrt2}&\frac{1}{\sqrt2}&\frac{1}{\sqrt2}&\frac{1}{\sqrt2}\\
\frac{1}{\sqrt2}&\frac{1}{\sqrt2}&0&0&
\frac{1}{\sqrt2}&\frac{1}{\sqrt2}&\frac{1}{\sqrt2}&\frac{1}{\sqrt2}\\
\frac{1}{\sqrt2}&\frac{1}{\sqrt2}&0&0&
\frac{1}{\sqrt2}&\frac{1}{\sqrt2}&\frac{1}{\sqrt2}&\frac{1}{\sqrt2}\\
\frac{1}{\sqrt2}&\frac{1}{\sqrt2}&\frac{1}{\sqrt2}&\frac{1}{\sqrt2}&
0&0&\frac{1}{\sqrt2}&\frac{1}{\sqrt2}\\
\frac{1}{\sqrt2}&\frac{1}{\sqrt2}&\frac{1}{\sqrt2}&\frac{1}{\sqrt2}&
0&0&\frac{1}{\sqrt2}&\frac{1}{\sqrt2}\\
\frac{1}{\sqrt2}&\frac{1}{\sqrt2}&\frac{1}{\sqrt2}&\frac{1}{\sqrt2}&
\frac{1}{\sqrt2}&\frac{1}{\sqrt2}&0&0\\
\frac{1}{\sqrt2}&\frac{1}{\sqrt2}&\frac{1}{\sqrt2}&\frac{1}{\sqrt2}&
\frac{1}{\sqrt2}&\frac{1}{\sqrt2}&0&0
\end{pmatrix}
=297
\end{equation*}
quasiequilibria.  
\begin{figure}
\label{cyclefour}
\begin{center}
\includegraphics{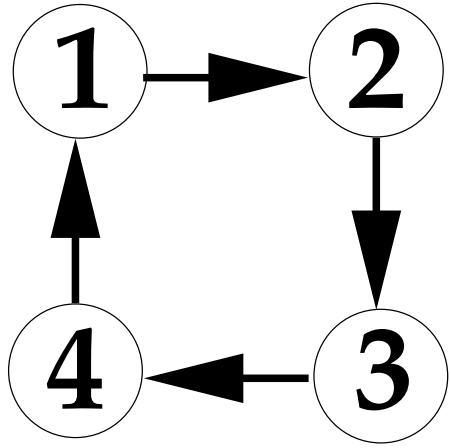}
\end{center}
\caption{Graphical game}
\end{figure}

\begin{figure}
\label{doublefourcycle}
\begin{center}
\includegraphics{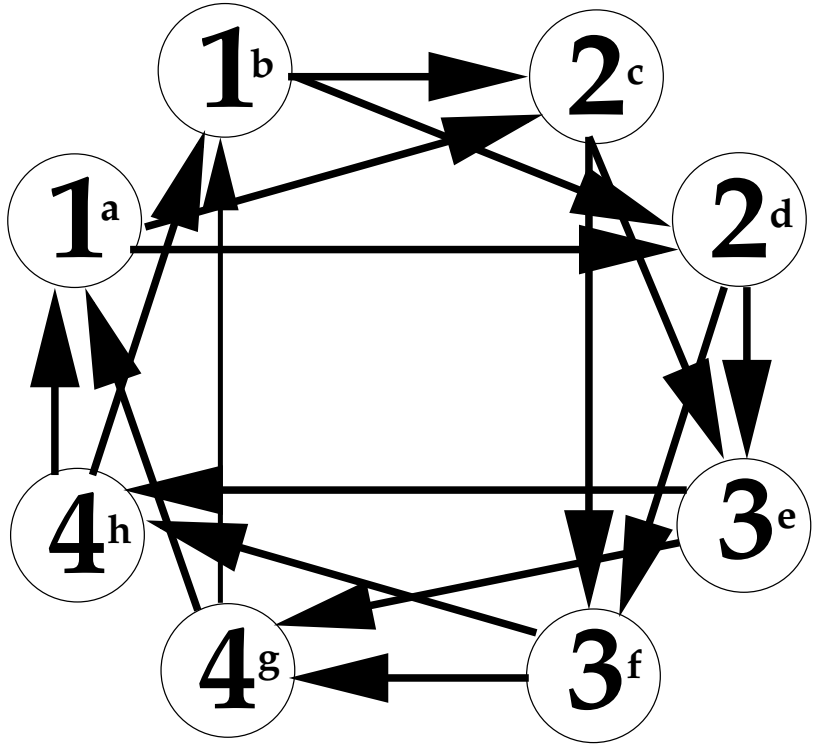}
\end{center}
\caption{Associated polynomial graph for graphical game}
\end{figure}

But suppose now that game obeys a graphical model as in
Figure 3.1. 
The nodes in the graphical model refer to the
players, and the edges specify that the payoff to the source player depends
on the actions of the target player.  For brevity, write
$a=\sigma_1(s_{11})$,
$b=\sigma_2(s_{12})$, $c=\sigma_2(s_{21})$, $d=\sigma_2(s_{22})$,
$e=\sigma_3(s_{31})$, $f=\sigma_3(s_{32})$, $g=\sigma_4(s_{41})$, and
$h=\sigma_4(s_{42})$.  Since the payoff to
player $1$ depends only on the actions of player $2$, equating the payoff
to player $1$ from pure strategies $s_{10}$ and $s_{11}$ gives
\begin{eqnarray*}
\lefteqn{u_1(s_{10},s_{20},\bullet)\sigma_2(s_{20})+u_1(s_{10},s_{21},\bullet)\sigma_2(s_{21})
+u_1(s_{10},s_{22},\bullet)\sigma_2(s_{22})}\\
&=&u_1(s_{11},s_{20},\bullet)\sigma_2(s_{20})+u_1(s_{11},s_{21},\bullet)\sigma_2(s_{21})
+u_1(s_{11},s_{22},\bullet)\sigma_2(s_{22})
\end{eqnarray*}
or
\begin{eqnarray*}
\lefteqn{\left(u_1(s_{11},s_{20},\bullet)-u_1(s_{10},s_{20},\bullet)\right)(1-c-d)+}\\
&&+\left(u_1(s_{11},s_{21},\bullet)-u_1(s_{10},s_{21},\bullet)\right)c
+\left(u_1(s_{11},s_{22},\bullet)-u_1(s_{10},s_{22},\bullet)\right)d=0.
\end{eqnarray*}
Thus for player $1$ we have two equations of the form
$$\bullet c+\bullet d+\bullet=0,$$
for player $2$ we have two equations of the form
$$\bullet e+\bullet f+\bullet=0,$$
for player $3$ we have two equations of the form
$$\bullet g+\bullet h+\bullet=0,$$
and for player $4$ we have two equations of the form
$$\bullet a+\bullet b+\bullet=0.$$
Then the associated polynomial graph is depicted in Figure \ref{doublefourcycle}.
The equation associated with the node labelled $1a$ equates the payoffs to
player $1$ from choosing $s_{11}$ (which $1$ does with probability $a$) or
choosing $s_{10}$.
The game has
\begin{equation*}
\hbox{per}
\begin{pmatrix}
0&0&\frac{1}{\sqrt2}&\frac{1}{\sqrt2}&
0&0&0&0\\
0&0&\frac{1}{\sqrt2}&\frac{1}{\sqrt2}&
0&0&0&0\\
0&0&0&0&
\frac{1}{\sqrt2}&\frac{1}{\sqrt2}&0&0\\
0&0&0&0&
\frac{1}{\sqrt2}&\frac{1}{\sqrt2}&0&0\\
0&0&0&0&
0&0&\frac{1}{\sqrt2}&\frac{1}{\sqrt2}\\
0&0&0&0&
0&0&\frac{1}{\sqrt2}&\frac{1}{\sqrt2}\\
\frac{1}{\sqrt2}&\frac{1}{\sqrt2}&0&0&
0&0&0&0\\
\frac{1}{\sqrt2}&\frac{1}{\sqrt2}&0&0&
0&0&0&0
\end{pmatrix}
=1
\end{equation*}
quasiequilibrium.  Indeed, this will always be the case for a graphical
model which is a directed cycle, where each player has the same number of
pure strategies.  The reason is that the indifference equations in this
case are {\em linear}, as we saw in this example.

The polynomial graph $G$ as defined in Theorem \ref{nashbernstein} contains
more refined information than the graphical model.  The partition into the
$T_i$'s also can be more refined than the partition of the set of all pure
strategies into the sets of pure strategies for each player.  Next we will
see an example of such a refinement when considering the reduction of
extensive-form games to normal-form, where actions correspond to branches
of the game tree.

\section{Extensive-Form Games}

Now we consider \emph{finite horizon extensive-form games}.  
(See for example \cite{OsborneRubinstein}, Chapter 6.)
Such a game takes place in a finite number of time steps, at each of which
only a single player (possibly nature) may move.  (Which player moves, and
what actions the player is allowed to take, may depend on what moves were
made previously.)  Such a game is completely specified as follows.  We
specify a set of players $I=\{1,\ldots,N\}$, and we specify a \emph{game
tree} $T$: a finitely branching tree of finite depth in which each non-leaf
node is labelled by a number in $0,\ldots,N$, each leaf is labelled by an
$N$-tuple of real numbers, and each branch emanating from a (non-leaf) node
labelled $0$ is assigned a positive real weight, so that the total
weight emanating from such a node is $1$.  (We consider the branches of this
tree to be directed away from the root.)  

Game play proceeds as follows. Each node of the tree represents a state of
the environment.  At each time step, if we are at a non-leaf node labelled
by $i$ in $1,\ldots,N$, then player $i$ acts by choosing one of the
branches emanating from that node.  Then the environment undergoes the
transition to the node at the end of that branch, and we advance to the
next time step.  If we are at a non-leaf node labelled by $0$, then the
environment instead makes a random transition along one of the branches
emanating from that node; the probability of each branch is given by its
weight.  If we are at a leaf node $\lambda$ labelled by
$(u_1(\lambda),\ldots,u_N(\lambda))$, then the
game is over, and each player $i$ accrues utility $u_i(\lambda)$.  Thus,
the leaf nodes are the outcomes of the game.

Unless otherwise specified, we will assume perfect information.  Not only
do all players have common knowledge of the specification of the game, but
whenever a player is about to move, that player knows what moves have been
made by all the other players (including nature) up to that point.

Every extensive-form game is equivalent to a normal-form game.  For each
node $\nu$ of the game tree, we write $E(\nu)$ for the set of edges
emanating from $\nu$.  Then the set of pure strategies of player $i$ is
$$S_i=\prod_{\nu\in T\atop{\rm label}(\nu)=i}E(\nu).$$
Given a pure strategy profile $s\in S=\prod_{i\in I}S_i$, we can compute
the probability of each leaf node $\lambda$ of the game tree.  A unique
path $\nu_0\nu_1\ldots\nu_m=\lambda$ leads from the root $\nu_0$ of $T$ to
$\lambda$.  Then $\Pr[\lambda|s]=\prod_{j=0}^{m-1}\Pr[\nu_j\to\nu_{j+1}|s]$,
where
$$\Pr[\nu_j\to\nu_{j+1}|s]=
\begin{cases}
1,&\hbox{$\nu_j$ is labelled by $i\in I$ and ${s_i}_{\nu_j}=(\nu_j\to\nu_{j+1})$}\\
0,&\hbox{$\nu_j$ is labelled by $i\in I$ and ${s_i}_{\nu_j}\neq(\nu_j\to\nu_{j+1})$}\\
{\rm wt}(\nu_j\to\nu_{j+1}),&\hbox{$\nu_j$ is labelled by $0$}\cr
\end{cases}$$
and so the utility functions of the normal-form game are given by
$$u_i(s)=\sum_{\lambda\in T\atop\lambda{\rm leaf}}u_i(\lambda)\Pr[\lambda|s].$$

We note that the game specification implies certain equalities among the
numbers $u_i(s)$.  
If we consider the set of normal-form games with a
fixed set of players $I$ and outcomes $S$ to be a linear space with basis
$\{u_i(s)\colon s\in S\}$, then the extensive-form games with the
same set of players $I$ and a fixed game tree having $S$ as the set of
outcomes lie in a linear subspace of this space, given by these
equalities.
Let $A$ be the set of non-leaf nodes of the tree which
are not labelled by $0$.  Then we can identify $S$ with $\prod_{\nu\in
A}E(\nu)$.  For any $s\in S$ and $\nu\in A$, we write $s_\nu={s_i}_\nu$,
where $i$ is the label of $\nu$.  Suppose $\nu\in A$ is an ancestor of
$\mu\in A$.  Then $\nu$ has a unique child $\alpha$ that is also an
ancestor of $\mu$ (possibly $\mu$ itself).  Let $\beta$ be any other child
of $\nu$.  If $s,s'\in A$ with $s_\nu=(\nu\to\beta)$ and
$$s'_\nu=
\begin{cases}
e,&\nu=\mu\cr
s_\nu,&\hbox{otherwise}\cr
\end{cases}$$
for some edge $e\in E(\mu)$, then $u_i(s)=u_i(s')$.  This is because
$\Pr[\lambda|s]=\Pr[\lambda|s']=0$ unless $\lambda$ is a descendant of
$\beta$ or $\lambda$ is not a descendant of $\nu$, and in either case
$\lambda$ cannot be a descendant of $\mu$.  In short, the node $\mu$ is
never reached, so it doesn't matter which action is chosen there.

If different players act at $\nu$ and $\mu$, then there is no way to
eliminate this redundancy, but when the same player $i$ acts at $\nu$ and
$\mu$, we can do so.  In this case we replace all the pure strategies which
are forced to be equal by a single pure strategy, called a \emph{reduced
pure strategy}.  See for example \cite{OsborneRubinstein}, p.~94.

We note that after iterated elimination of strictly dominated pure
strategies, for any node all of whose children are leaves, the payoffs to
the player who acts at that node must be equal at all these child leaves.
If nature acts at such a node $\nu$ whose children are leaves
$\lambda_1,\ldots,\lambda_k$, then we can replace $\nu$ by a leaf with
utilities $u_i(\nu)=\sum_{l=1}^k \hbox{wt}(\nu\to\lambda_l)u_i(\lambda)$
for each $i\in I$.  So we assume nature never acts at such nodes.

For extensive-form games, the equilibrium concept can be refined.  Each
subtree of the game tree induces a new extensive-form game, called a
\emph{subgame}.  Each pure strategy of the original game induces a pure
strategy of each subgame by restriction to that subtree, and thus each
strategy profile of the original game induces a strategy profile of each
subgame.  A strategy profile is a \emph{subgame perfect Nash equilibrium}
of an extensive-form game if it induces a Nash equilibrium of each subgame.

We can find a subgame perfect pure strategy Nash equilibrium by
\emph{backwards induction}.  We construct the pure strategy profile as
follows.   We perform iterated elimination of strictly dominated
strategies.  Then at each node all of whose children are leaves, we choose
one leaf (recall that the payoffs of all leaves for the player who acts at
that node will be the same).  We assign this branch to the corresponding
component of the pure strategy profile, replace this node by this leaf, and
repeat the procedure on the resulting subtree.

We begin our analysis of totally mixed Nash equilibria of extensive form games
by noting the following:

\begin{theorem}
\label{subgameperfect}
All totally mixed Nash equilibria of an extensive form game are subgame
perfect.
\end{theorem}
\begin{proof}
Let $\sigma$ be a totally mixed Nash equilibrium of an extensive form game
with $N$ players defined by game tree $T$. Note that the strategy profile
induced by $\sigma$ on every subgame is also totally mixed.  Let $\nu$ be a
non-leaf node of $T$.  Let $\tilde \sigma$ be the strategy profile induced
by $\sigma$ in the subgame induced by $\nu$.  Let $\tilde s_j$ and $\tilde
t_j$ be pure strategies of player $j$ in this subgame.  Choose an action
for $j$ at each node $\mu$ that is not a descendant of $\nu$ where $j$
acts, such that if $\mu$ is an ancestor of $\nu$ then $j$ chooses the
branch leading towards $\nu$, and use this choice to extend $\tilde s_j$
and $\tilde t_j$ to pure strategies $s_j$ and $t_j$ of player $j$ in the
original game.  (So, $s_j$ and $t_j$ specify the same actions outside the
subtree.)
Let $\nu_0\ldots\nu_m=\nu$ be the unique path from the root $\nu_0$
of $T$ to $\nu$.  We have $u_j(s_j,\sigma_{-j})=u_j(t_j,\sigma_{-j})$.  
Let $L$ be the set of all leaves of $T$ under $\nu$ and $L'$ be the set of
all other leaves.  Then 
\begin{eqnarray*}
u_j(s_j,\sigma_{-j})&=\sum_{\lambda\in
L}u_j(\lambda)\Pr[\lambda|(s_j,\sigma_{-j})]+\sum_{\lambda\in
L'}u_j(\lambda)\Pr[\lambda|(s_j,\sigma_{-j})]\\
&=\sum_{\lambda\in
L}u_j(\lambda)\Pr[\lambda|(s_j,\sigma_{-j})]+\sum_{\lambda\in
L'}u_j(\lambda)\Pr[\lambda|(t_j,\sigma_{-j})]
\end{eqnarray*} since $s_j$ and $t_j$
choose the same actions outside the subtree.  Thus 
\begin{equation}
\label{subtree}
\sum_{\lambda\in
L}u_j(\lambda)\Pr[\lambda|(s_j,\sigma_{-j})]=\sum_{\lambda\in
L}u_j(\lambda)\Pr[\lambda|(t_j,\sigma_{-j})].  
\end{equation}
Furthermore, for any
$\lambda\in L$, we have
\begin{eqnarray*}
\Pr[\lambda|(s_j,\sigma_{-j})]&=&
\Pr[\lambda|(\tilde s_j,\tilde\sigma_{-j})]
\prod_{k=0}^{m-1}\Pr[\nu_k\to\nu_{k+1}|(s_j,\sigma_{-j})]\\
&=&\Pr[\lambda|(\tilde s_j,\tilde\sigma_{-j})]
\prod_{k=0}^{m-1}\Pr[\nu_k\to\nu_{k+1}|(t_j,\sigma_{-j})].
\end{eqnarray*}
Noting that the common factor
$\prod_{k=0}^{m-1}\Pr[\nu_k\to\nu_{k+1}|(t_j,\sigma_{-j})]$ in equation
(\ref{subtree}) is positive by
our choice of $s_j,t_j$ and because $\sigma$ is totally mixed, we have
that 
\begin{eqnarray*}
u_j(\tilde s_j,\tilde\sigma_{-j})
&=&\sum_{\lambda\in L}u_j(\lambda)Pr[\lambda|(\tilde s_j,\tilde\sigma_{-j})]\\
&=&\sum_{\lambda\in L}u_j(\lambda)Pr[\lambda|(\tilde
t_j,\tilde\sigma_{-j})]\\
&=&u_j(\tilde t_j,\tilde\sigma_{-j}).  
\end{eqnarray*}
Thus $\tilde\sigma$ is a (totally mixed) Nash equilibrium of the subgame
induced by $\nu$.
\end{proof}

In light of this observation, the divide-and-conquer approach to finding
all Nash equilibria of a normal form game can be modified in the spirit of
backwards induction to finding all subgame perfect equilibria (including
mixed ones) of an extensive form game.  Recall that in a normal form game,
we would consider subproblems in which one pure strategy of one player $i$
was removed.  Now we instead consider subproblems in which, for some edge
$\nu\to\mu$ where $i$ acts at $\nu$, we delete that edge and the entire
subtree below $\mu$.  We compute the normal form for the game
described by this pruned tree and recursively find all its subgame perfect
equilibria.  Each such equilibrium $\sigma$ induces an equilibrium
$\tilde\sigma$ in the subgame under $\nu$ in the pruned tree.  To check
whether $\sigma$ is an equilibrium of the original game, we recursively
compute all the equilibria of the subgame under $\mu$ (where $i$ does not
act), and check that for each such equilibrium $\tau$, we have
$u_i(\tilde\sigma)\geq u_i(\tau)$.

We saw during the above proof that for a totally mixed strategy profile
$\sigma$, the equations $u_j(s_j,\sigma_{-j})=u_j(t_j,\sigma_{-j})$ for all
pure strategies $s_j,t_j$ of $j$ imply the corresponding equations for each
subtree.  The converse implication also clearly holds.  

We will now associate a polynomial graph to a system of equations for the
quasiequilibria of an extensive-form game, so that we can apply Theorem
\ref{nashbernstein}.  For each node in the game tree where a player acts,
we will have a variable for every edge emanating from that node except one
distinguished edge.  This is
because the sum of the probabilities of choosing each of those edges must
be $1$, so we eliminate one variable.  Thus, we compare the payoffs between
choosing the distinguished edge and choosing any other edge.  The equations
will be indifference equations for subgames of the extensive-form game.

\begin{theorem}
\label{extensivetheorem}
The set of quasiequilibria of a generic extensive-form game is either empty or has
positive dimension.
\end{theorem}
\begin{proof}
Consider an extensive form game with players $I=1,\ldots,N$ and game tree
$T$.  Let $A$ be the set of non-leaf nodes in $T$ not labelled by $0$.
For each $\nu\in A$, let $E(\nu)$ be
the set of edges emanating from $\nu$.  For each $\nu\in A$, let $i$ be the
player which acts at $\nu$ and pick an element $e_{i\nu}\in E(\nu)$.
Let $d=\sum_{\nu\in A}|E(\nu)-1|$  and
partition $d$ as $\coprod_{\nu\in A}\left(E(\nu)-\{e_{i\nu}\}\right)$.  
Define a directed graph $G$ on a set of
$d$ vertices $$\bigcup_{\nu\in A}\left\{n_e\mid
e\in E(\nu)-\left\{e_{i\nu}\right\}\right\}$$ as
follows: there is an edge from $n_e$ with $e\in
E(\nu)-\{e_{i\nu}\}$ to
$n_{e'}$ with $e'\in E(\mu)-\{e_{j\mu}\}$ if $i\neq j$,
$\nu$ is an ancestor of $\mu$, either $e$ or $e_{i\nu}$ lies on the path
from $\nu$ to $\mu$, and if $i$ acts at some node $\kappa$
between $\nu$ and $\mu$, then the edge $e_{i\kappa}$ lies on the path from
$\nu$ to $\mu$.  We will define a system of equations equivalent to the
equations defining totally mixed Nash equilibria
of the extensive form game and satisfying conditions 1 to 3 of Theorem
\ref{nashbernstein}.  The polynomial graph $G$ is acyclic, so Corollary
\ref{cyclicgraph} implies our assertion.

First we must state what the equations are.  Fix a node $\nu\in A$ and let
$i$ be the player which acts at $\nu$.  Then
$|E(\nu)|-1$ equations refer to the subgame induced by this node.  
For each $e\in E(\nu)$, define the pure strategy $s_{ie}$ of $i$ in this
subgame by $s_{ie}(\nu)=e$ and $s_{ie}(\mu)=e_{i\mu}$ for any node $\mu$
below $\nu$ where $i$ acts.  
Writing $\tilde\sigma$ for the strategy profile induced by $\sigma$ in the
subgame under $\nu$, the $|E(\nu)|-1$ equations are the equations 
$u_i(s_e,\tilde\sigma_{-i})=u_i(s_{e_{i\nu}},\tilde\sigma_{-i})$ for $e\in
E(\nu)-\{e_{i\nu}\}$.  In these equations we eliminated $\sigma(e_{j\mu})$
for every $\mu$ below $\nu$ where $i$ does not act, by substituting 
$1-\sum_{e\in\tilde E(\mu)-e_{\{j\mu}\}}\sigma_j(e)$ for
$\sigma_j(e_{j\mu})$.

These are some of the indifference equations for the subtree below $\nu$,
which as we saw in the previous theorem are implied by the
indifference equations for the whole tree.  We show by induction that these
equations also imply all the indifference equations for the subtree below
$\nu$.  (Thus we will have the indifference equations for every subtree,
and hence the whole tree, i.e., the original game.) Firstly, $i$ is
indifferent between {\em all} $i$'s pure strategies in the subgame below
$\nu$, because although we fixed $i$'s pure strategies at nodes $\mu$ below
$\nu$ where $i$ acts to be $e_{i\mu}$, we also have that $i$ is indifferent
between $i$'s pure strategies in the subgame below $\mu$ by the induction
hypothesis.  Secondly, consider any other player $j$. Let
$\mu_1,\ldots,\mu_m$ be the nodes below $\nu$ where $j$ acts, such that $j$
does not act at any node between $\nu$ and $\mu_k$ for any $k$.  Let
$\tilde s_j,\tilde t_j$ be pure strategies of $j$ in the subgame below
$\nu$, and write $\tilde {s_j}_k,\tilde {t_j}_k$ for the respective induced
pure strategies of $j$ in the subgame below $\mu_k$.  So $\tilde
s_j=(\tilde {s_j}_1,\ldots,\tilde {s_j}_m)$ and $\tilde
t_j=(\tilde{t_j}_1,\ldots,\tilde{t_j}_m)$.  Write the set $L$ of leaves below
$\nu$ as $L=L_0\cup\bigcup_{k=1}^m L_k$, where $L_0$ is the set of leaves
$\lambda$ such that $j$ does not act between $\nu$ and $\lambda$ and $L_k$
is the set of leaves below $\mu_k$ for $k=1,\ldots,m$.  Then 
\begin{eqnarray*}
u_j(\tilde s_j,\tilde \sigma_{-j})&=&\sum_{\lambda\in
L}u_j(\lambda|\tilde s_j,\tilde\sigma_{-j})\\
&=&\sum_{\lambda\in
L_0}u_j(\lambda|\tilde\sigma_{-j})+\sum_{k=1}^m\sum_{\lambda\in
L_k}u_j(\lambda|\tilde {s_j}_k,\tilde\sigma_{-j})\\
&=&\sum_{\lambda\in
L_0}u_j(\lambda|\tilde\sigma_{-j})+\sum_{k=1}^m\sum_{\lambda\in
L_k}u_j(\lambda|\tilde {t_j}_k,\tilde\sigma_{-j})\\
&=&u_j(\tilde t_j,\tilde\sigma_{-j})
\end{eqnarray*}
since for each $k$, 
$\sum_{\lambda\in L_k}u_j(\lambda|\tilde {s_j}_k,\tilde\sigma_{-j})=
\sum_{\lambda\in L_k}u_j(\lambda|\tilde {t_j}_k,\tilde\sigma_{-j})$
by the induction hypothesis.

We can already see that the set of solutions to these equations, if
nonempty, is positive-dimensional.  If player $i$ acts at the root $\nu$,
then for any edge $e$ emerging from $\nu$, $\sigma_i(e)$ does not appear in
any of the equations.

All the monomials occurring in these equations are squarefree.  For each
leaf $\lambda$ under $\nu$, let the path from $\nu$ to $\lambda$ be
$\nu=\nu_1\ldots\nu_k=\lambda$.  Then for any player $j$ with pure strategy
$\tilde s_j$, we have $\Pr[\lambda|\tilde
s_j,\tilde\sigma_{-j}]=\prod_{l=1}^{k-1}\Pr[\nu_l\to\nu_{l+1}|\tilde
s_j,\tilde\sigma_{-j}]$, and each nonconstant term in the product is
$\sigma_n(\nu_l\to\nu_{l+1})$ for some player $n\neq j$.  So for any edge $e$
where $n$ acts, the
variable $\sigma_n(e)$ occurs at most once in such a product.  In fact
$\sigma_n(e)$ occurs in such a product for at most one $e\in E(\nu_l)$.
(That is, if $e,e'\in E(\nu_l)$ then $\sigma_n(e)$ and $\sigma_n(e')$ do
not both occur in this monomial.  So condition 2 of Theorem
\ref{nashbernstein} holds.)
When we eliminate $\sigma_n(e_{n\nu_l})$, we replace it by an affine
expression, so this remains true.
Thus condition 1 of Theorem \ref{nashbernstein} holds.

The equations corresponding to $E(\nu)-\{e_{i\nu}\}$ concern only the
subgame below $\nu$, so $\sigma_j(\mu\to\kappa)$ occurs in these equations
only if $\nu$ is an ancestor of $\mu$.  Furthermore, if $i$ acts at $\kappa$
below $\nu$, then $\sigma_i(e)$ does not occur for any edge $e\in
E(\kappa)-\{e_{i\kappa}\}$, since we fix that $i$ chooses $e_{i\kappa}$.
For the same reason $\sigma_j(e)$ does not occur for $e\in
E(\mu)-\{e_{j\mu}\}$ for any $\mu$ that lies below $\kappa$ but not below
$e_{i\kappa}$.  Thus condition 3 holds.  
\end{proof}

Our result does not contradict Harsanyi's generic finiteness theorem
\cite{harsanyi}, because generically, iterated elimination of weakly
dominated strategies/backward induction will lead to a unique subgame
perfect equilibrium (and so indeed there will be no totally mixed Nash
equilibria).  On the other hand, another way to look at our result is that
in every {\em interesting} extensive-form game---one which is not
completely solved by backward induction, giving a unique equilibrium---the
set of totally mixed Nash equilibria is also interesting; it has positive
dimension.

In particular, if $\nu$ is a node all of whose children are leaves, the
equations corresponding to $\nu$ will be equations between constants,
stating that for the player $i$ who acts at $\nu$, the utilities
$u_i(\lambda)$ at all the leaves $\lambda$ below $\nu$ must be equal.  This
is true if iterated elimination of strictly dominated pure strategies has
already been performed on this game.  

It is clear that the system of equations we obtained is not canonical,
since we have made arbitrary choices of the edges $e_{i\nu}$ and the
subtrees below each possible choice are different.  Choosing a different
system may make it easier to compute the set of quasiequilibria.

\begin{figure}
\begin{center}
\includegraphics{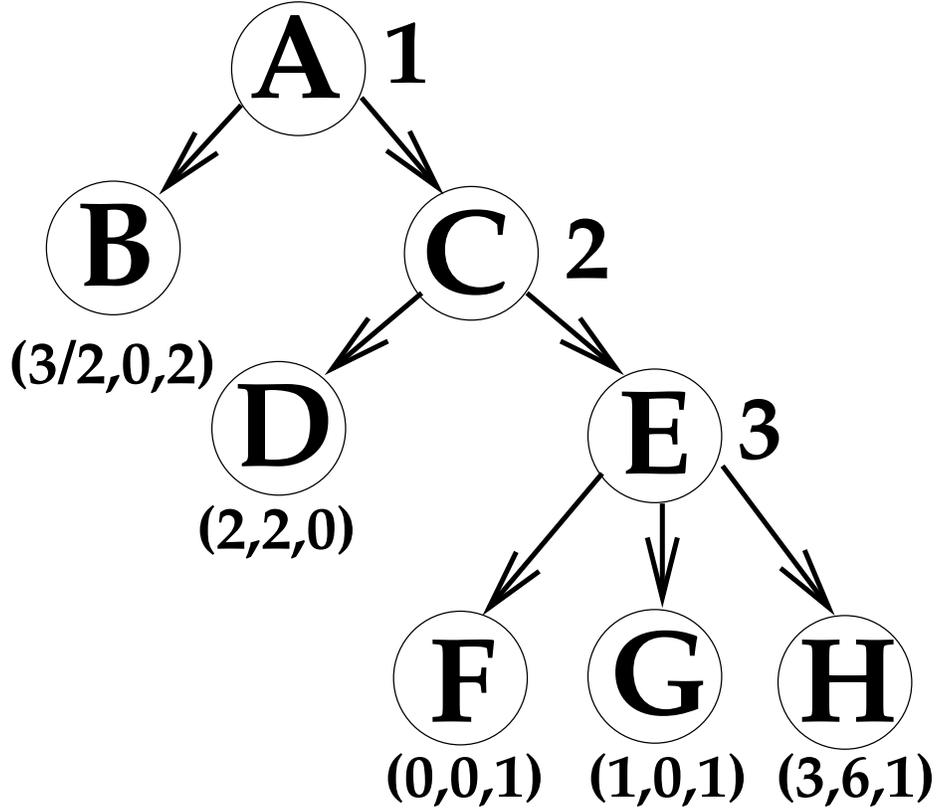}
\end{center}
\label{hyperbolic}
\caption{An Extensive Form Game}
\end{figure}

We now present an example where the set of totally mixed Nash equilibria is a
positive-dimensional semialgebraic variety.  Consider the extensive form
game specified in Figure \ref{hyperbolic}.
The polynomial graph associated with this game tree is depicted in Figure
3.4. 
\begin{figure}
\begin{center}
\includegraphics{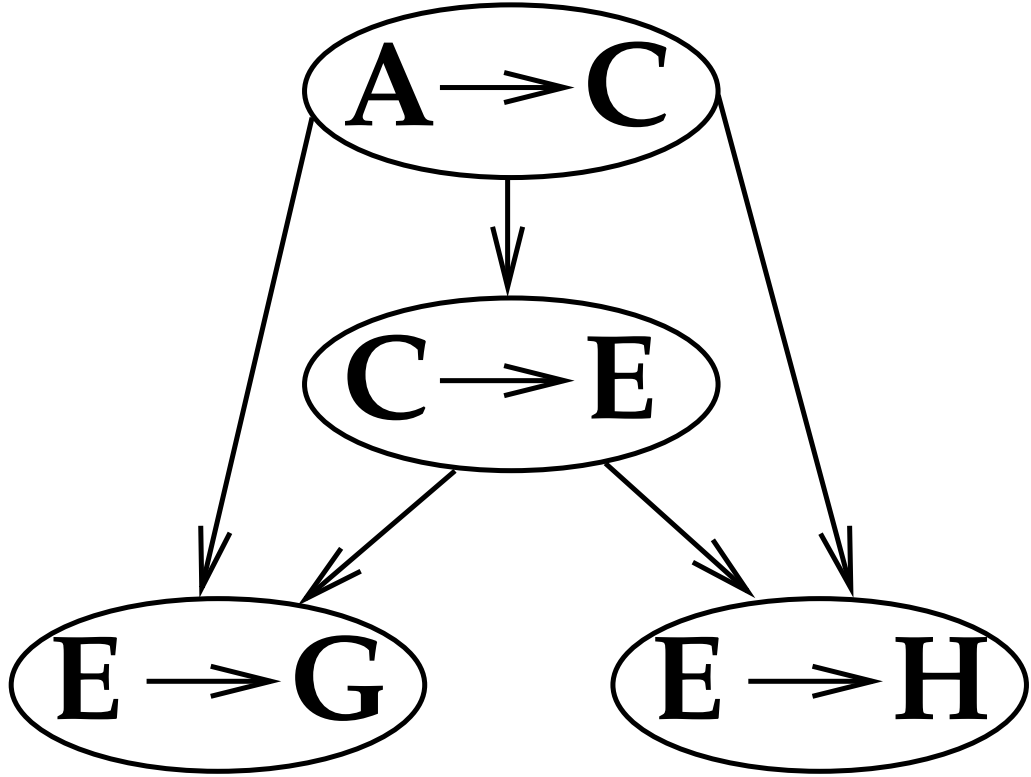}
\end{center}
\label{hyperbolicgraph}
\caption{Associated Polynomial Graph For An Extensive Form Game}
\end{figure}
For brevity, we write for example $\sigma_1(C)$ for $\sigma_1(A\to C)$.
The quasiequilibria obey a system of 4 equations as in Theorem
\ref{extensivetheorem}.  The equation associated with the edge $E\to G$
equates the payoff to player $3$ from choosing this edge with that from
choosing the edge $E\to F$, i.e., $u_3(F)=u_3(G)$.  No variables occur in
this equation, that is, it is an equation between constants.  Similarly,
the equation associated with the edge $E\to H$ is $u_3(F)=u_3(H)$.  The
equation associated with the edge $C\to E$ is $u_2(D)=u_2(E)$, where we
have written $u_2(E)$ for the expected payoff $u_2(E,\sigma_{-2})$ to
player 2 for choosing the edge $C\to E$, given the strategy profile of the
other players.  In this case
$u_2(E)=u_2(F)\sigma_3(F)+u_2(G)\sigma_3(G)+u_2(H)\sigma_3(H)$, so
$$u_2(D)=u_2(F)+\left(u_2(G)-u_2(F)\right)\sigma_3(G)+\left(u_2(H)-u_2(F)\right)\sigma_3(H).$$
Finally, the equation associated to the edge $A\to C$ is
\begin{eqnarray*}
u_1(B)&=&u_1(C)\\
&=&u_1(D)\left(1-\sigma_2(E)\right)+u_1(F)\sigma_2(E)\left(1-\sigma_3(G)-\sigma_3(H)\right)\\
&&\quad+u_1(G)\sigma_2(E)\sigma_3(G)+u_1(H)\sigma_2(E)\sigma_3(H).
\end{eqnarray*}

Looking at the specific payoffs in Figure \ref{hyperbolic},
we see that the payoffs to player 3 for choosing $F$, $G$, or $H$ are equal, as
required.  Equating the payoffs to player 2 for choosing $D$ or $E$, we get
$6\sigma_3(H)=2$, or $\sigma_3(H)=\frac{1}{3}$.  This leaves $\sigma_3(G)$
free to vary such that $0<\sigma_3(G)<\frac{2}{3}$.  Finally, we must
equate the payoffs to player 1 for choosing $B$ or $C$.  This gives
$$2(1-\sigma_2(E))+\sigma_2(E)\left(\sigma_3(G)+1\right)={3\over2}$$
or
$$\sigma_2(E)(1-\sigma_3(G))={1\over 2}.$$
Thus the points $\sigma_3(G)$ and $\sigma_2(E)$ lie on a hyperbola.  This
hyperbola intersects the interior of the product of simplices.  For
instance, the point $\sigma_3(G)=\frac{5}{12}$ (so $\sigma_3(F)=\frac{1}{4}$)
and $\sigma_2(E)=\frac{6}{7}$ lies in this intersection.  So the set of
quasiequilibria is a portion of a hyperbolic cylinder, the product of a
segment of a hyperbola with a line segment (since $\sigma_1(B)$ varies
freely with $0<\sigma_1(B)<1$).

We can analyze this game a little further.  Player 3 would like player 1 to
sometimes choose $B$, but cannot force player 1 always to choose $B$, since
if player 2 always chooses $D$ then both player 1 and player 2 are better
off with player 1 choosing $C$.  The best player 3 can do is make the
payoffs to player 1 from choosing $B$ and $C$ equal.  Now if player 3 made
player 2 get a greater payoff from choosing $D$ than $E$, then player 2
would always choose $D$, player 1 would always choose $C$, and player 3
would get nothing.  So player 3 must make $u_2(D)\leq u_2(E)$.  We analyzed
the case $u_2(D)=u_2(E)$ above.  If player 3 makes
$\sigma_3(H)>\frac{1}{3}$, then $u_2(D)<u_2(E)$ and player 2 will always
choose $E$.  Then the payoff to player 1 from choosing $C$ is
$\sigma_3(G)+3\sigma_3(H)$.  Thus we have
$\sigma_3(G)+3\sigma_3(H)=\frac{3}{2}$ with
$\frac{1}{3}<\sigma_3(H)\leq\frac{1}{2}$ (this makes
$0\leq\sigma_3(G)<\frac{1}{2}$ and
$\frac{1}{6}<\sigma_3(F)\leq\frac{1}{2}$).  Then $\sigma_1(C)$ varies
freely with $0\leq\sigma_1(C)\leq1$, so we have a rectangle of partially
mixed equilibria.  Player 3 is better off choosing these, since then the
outcome $D$ where player 3 gets zero payoff is never reached.  Along the
line $\sigma_3(G)+3\sigma_3(H)=\frac{3}{2}$, equilibria with greater
$\sigma_3(H)$ \emph{Pareto dominate} those with smaller $\sigma_3(H)$,
i.e., they make some player better off and no player worse off.
Specifically, the payoff to player 2 increases, the payoff to player 1 is
always $\frac{3}{2}$, and the payoff to player 3 stays the same at
$2(1-\sigma_1(C))+\sigma_1(C)=2-\sigma_1(C)$.  Thus the Pareto dominant
equilibrium among those on this line is that player 3 has
$\sigma_3(F)=\frac{1}{2}$, $\sigma_3(G)=0$, and $\sigma_3(H)=\frac{1}{2}$.
On the other hand, at the pure strategy equilibrium where player 3 always
chooses $H$, we have that player 1 always chooses $C$, and the payoff to
player 3 falls from $2-\sigma_1(C)$ to $1$.  Thus player 3 does not prefer
this equilibrium, and instead mixes $F$ and $H$ equally to have some chance
of a higher payoff.  As $\sigma_1(C)$ increases, the payoff to player 3
decreases and the payoff to player 2 increases, so the equilibria along
this line do not Pareto dominate each other.  Thus without introducing
other issues (such as risk-aversion) there is no criterion for predicting
which of the equilibria along the line $0<\sigma_1(C)<1$,
$\sigma_2(E)=1$, $\sigma_3(F)=\sigma_3(H)=\frac{1}{2}$ should be chosen.

%

\section{Games With Emergent Node Tree Structure}

So far we have been discussing normal form games with finite numbers of
players, each with a finite number of pure strategies.  Such a game is
defined by giving a set of players $I=\{1,\ldots,\N\}$, for each player $i$
a finite set of pure strategies $S_i$, and for each pure strategy profile
$\sigma$ (element of the product $S=\prod_{i\in I}S_i$) and each player $i$
the utility $u_i(\sigma)$ received by that player when that strategy
profile is played.  Now we will introduce a particular kind of structure that
a normal form game may have.  

We now define an \emph{emergent node tree structure} on a normal form game.
This is a new model for games in which the players can be hierarchically
decomposed into groups.  Usually such hierarchical decomposition is
discussed in the framework of cooperative game theory.  Instead, we define
certain conditions on the payoff functions in a noncooperative game such
that a given hierarchical decomposition ``makes sense'', in a way that we
will define precisely.  At the end of this section we briefly describe how
our framework relates to that of cooperative game theory.

{\bf Definition.} An \emph{emergent node tree structure} on a normal form game with player
$I=\{1,\ldots,\N\}$, pure strategy sets $S_i$ for $i\in I$, and utility
functions $u_i\colon\prod_{i\in I}S_i\to\R$ to consist of:

\begin{itemize}
\item A tree $T$ with $N$ leaves.  The leaves are in bijection with the players
$I=\{1,\ldots,\N\}$.  Write $C_v$ for the set of children of a node $v\in
T$, $B_v$ for the set of its siblings, and $f(v)$ for its parent.

\item For each non-leaf, non-root node $v$ of the tree (which we call an
{\em emergent player}), a set $S_v$ of pure strategies, with $|S_v|\leq
\prod_{w\in C_v}|S_w|$.

\item For each non-leaf, non-root node $v$, for each element $s_{C_v}$ of the
product $S_{C_v}=\prod_{w\in C_v} S_w$ of the pure strategies of its children and
each element $s_{vk}$ of $S_v$, a number $p_v(k,s_{C_v})$
signifying the probability that the (emergent) strategy of the emergent
player $v$ is $s_{vk}$ when the strategies of its children are given by
$s_{C_v}$.  So if $v$ has $K$ pure strategies, then $\sum_{k=1}^K
p_v(k,s_{C_v})=1$.  If the children of $v$ execute a mixed
strategy, then the emergent mixed strategy of $v$ is given by
multilinearity.  Thus we have defined a linear map from the strategy space
of the children to the strategy space of the parent.  We require that this
map have full rank.

\item For each non-root node $v$ (including the leaf nodes), real numbers
$\gamma_{vw}$ for each non-root ancestor $w$ of $v$ and real numbers
$U_v(s)$ for each element $s\in S_v\times \prod_{w\in B_v}S_w$.  From these
we define a utility function $u_v$, which is a sum of two terms:
$U_v(\sigma_{v,B_v})$, a multilinear function of the strategies executed by
$v$ and its siblings in $B_v$, and $\sum_{{\rm
nonroot\,\,ancestors}\,\,u}\gamma_{vw} u_w$.  We require that the utility
function $u_v$ at a leaf node $v$ be equal to the utility function $u_i$ of
the player $i$ corresponding to the leaf node $v$.  \end{itemize}

We will refer to an emergent node tree structure as an \emph{ENT} for
short.  Note that for a given normal form game, we can always define a
class of ENTs by defining a tree with a single emergent node (the root
node), so that all the leaf nodes are siblings.  We call such an ENT
\emph{trivial}.  For any given normal form game, there need be no
nontrivial ENT, or there may be many distinct possible ENTs.

The behavior of the \emph{emergent players} is completely determined by the
behavior of the actual players (the leaf nodes).  The \emph{emergent
strategy} $\sigma_v$ executed by the emergent player $v$ when the actual
players execute strategy profile $\sigma$ is defined recursively by
multilinearity:
$$\sigma_v(s_{vk})=\sum_{s\in S_{C_v}}p_v(k,s)\prod_{w\in
C_v}\sigma_w(s_w).$$
So we compute the emergent strategies from the bottom up.

From the above definition, we see that at a non-root node $w$ of the tree, the
utility function is
\begin{eqnarray*}
u_w(\sigma)&=&U_w(\sigma_w)+\sum_{{\rm
nonroot\,\,ancestor}\,\,v}\gamma_{wv}u_{v}(\sigma)\\
&=&\sum_{s\in
S_w\times\prod_{x\in B_w}S_x}U_w(s)\sigma_w(s_w)\prod_{x\in
B_w}\sigma_x(s_x)+\sum_{{\rm nonroot\,\,ancestor}\,\,v}\gamma_{wv}u_{v}(\sigma).
\end{eqnarray*}
So we compute the utility from the top down.

We see that the utilities of each actual player (the leaf nodes) may depend
on the strategies executed by every other actual player.  So, the graphical
model of the actual game may be the complete graph.  Imposing an emergent
node tree structure, corresponds to  deleting some of these edges and
adding more nodes, and edges connected to those nodes, to the graph, so
that the new graph has a nontrivial structure.  With the addition of the
new variables $\sigma_v(s_{vk})$, we get more information about the
sparsity of our multilinear equations.  

In our definition, we did not require that the numbers $\gamma_{vw}$ have the
same sign for all descendants $v$ of a node $w$.  Thus, our definition does
not require that the emergence of a node represent a common interest among
its descendant nodes (although of course it does cover that situation).  

\begin{figure}
\begin{center}
\includegraphics{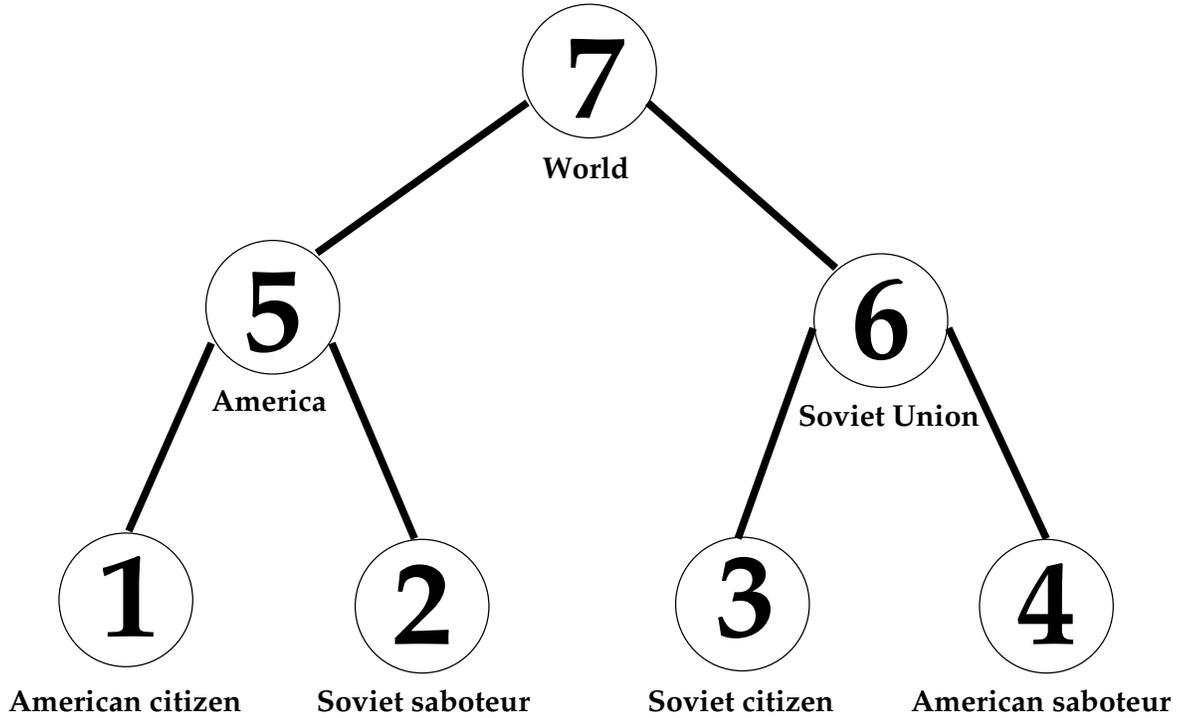}
\end{center}
\caption{Emergent Node Structure For The Saboteur Game}
\label{saboteur}
\end{figure}

For example, consider a normal form game with the ENT in Figure
\ref{saboteur} where the leaf nodes
correspond to
\begin{enumerate}
\item An American citizen
\item A Soviet saboteur living in America
\item A Soviet citizen
\item An American saboteur living in the USSR
\end{enumerate}
The parent of nodes $1$ and $2$ is node $5$, corresponding to America, the
parent of nodes $3$ and $4$ is node $6$, corresponding to the USSR, and the 
the root is node $7$, corresponding to the world.  Then while
$\gamma_{15}>0$ and $\gamma_{36}>0$, we have $\gamma_{25}<0$ and
$\gamma_{46}<0$.  

We now define a natural refinement of the equilibrium concept for games
with an ENTs.

{\bf Definition.}
If a normal form game has an ENT as defined above,
then a Nash equilibrium $\sigma$ of that game is \emph{hierarchically
perfect} with respect to this ENT if for every
emergent node $v$, given the strategies induced on the siblings of $v$ by
$\sigma$, the payoff $u(v)$ at $v$ cannot be increased by changing only
$\sigma(v)$.

Note that since our definition requires the linear map from the strategy
space of the children of $v$ to the strategy space of $v$ to be full-rank,
any strategy $\sigma'(v)$ deviating from $\sigma(v)$ which could result in
a higher payoff $u(v)$ would be achievable by some strategy profile of the
descendants of $v$.

We will also need the following definition:

{\bf Definition.}
A strategy profile of a normal form game with an ENT is \emph{totally
mixed with respect to this ENT} if it is totally mixed in the usual sense
and the emergent strategies at each emergent node are also totally mixed.

\begin{theorem}
For a generic game with an ENT as above, construct a
directed graphical model $G$ whose nodes are the nodes of the tree except
the root, with edges as
follows: the children in $T$ of a node $v$ form a directed clique in $G$,
and each such child also has a directed edge from $v$ and each ancestor of
$v$ except the root, and from each of their siblings.  Then the Bernstein number we obtain by
applying Theorem \ref{graphicalmodeltheorem} to this directed graphical model is an
upper bound on the number of totally mixed Nash
equilibria of this game which are hierarchically perfect and totally mixed
with respect to this ENT.
\end{theorem}
\begin{proof}
This is the graphical model we would obtain if all the emergent players
were actual players.  That is, we have ignored the
equations 
$$\sigma_v(s_{vk})=\sum_{s\in S_{C_v}}p_v(k,s)\prod_{w\in
C_v}\sigma_w(s_w).$$ 
So the set of totally mixed Nash equilibria of our game which are
hierarchically perfect with respect to this ENT is
a subset of the set of totally mixed Nash equilibria of the game with this
graphical model.
\end{proof}

Generically, there may be no hierarchically perfect totally mixed Nash
equilibria.  If the system of equations defining the quasiequilibria of the
game with the directed graphical model is $0$-dimensional, then none of the
finitely many solutions to this system may satisfy the
additional equations $$\sigma_v(s_{vk})=\sum_{s\in
S_{C_v}}p_v(k,s)\prod_{w\in C_v}\sigma_w(s_w).$$

\begin{figure}
\begin{center}
\includegraphics{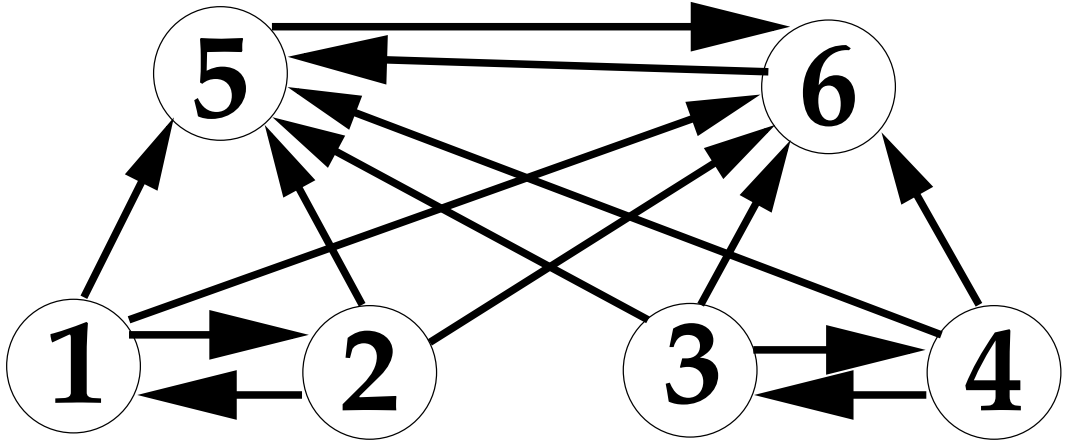}
\end{center}
\label{saboteurgraph}
\caption{Graphical Model For The Saboteur Game}
\end{figure}

For example, consider a game as in Figure \ref{saboteur} in which each
actual player has two pure strategies and each emergent player also has two
pure strategies.  Generically, a game with 4 players, each with 2 pure
strategies, would have 
\begin{equation*}
\hbox{per}
\begin{pmatrix}
0&1&1&1\\
1&0&1&1\\
1&1&0&1\\
1&1&1&0
\end{pmatrix}
=9
\end{equation*}
quasiequilibria.  On the other hand, if the game has an ENT as in Figure
\ref{saboteur}, then the directed graphical model given by the theorem is
as in Figure
\ref{saboteurgraph}.
Thus there is no more than
\begin{equation*}
\hbox{per}
\begin{pmatrix}
0&1&0&0&1&0\\
1&0&0&0&1&0\\
0&0&0&1&0&1\\
0&0&1&0&0&1\\
0&0&0&0&0&1\\
0&0&0&0&1&0\\
\end{pmatrix}
=1
\end{equation*}
quasiequilibrium which is hierarchically perfect and totally mixed with
respect to this ENT.  Indeed this would hold whenever the ENT is a binary
tree, that is, each non-leaf node has two children, and all siblings
have the same number of pure strategies. 

For example, say that if players 1 and 2 either both choose their 0th pure
strategy or both choose their 1st pure strategy, then the emergent strategy
of node 5 is $s_{51}$, otherwise it is $s_{50}$.  Similarly, if players 3
and 4 either both choose their 0th pure strategy or both choose their 1st
pure strategy, then the emergent strategy of node 6 is $s_{61}$, otherwise
it is $s_{60}$.  Let $U_5$ and $U_6$ be given by
\begin{equation}
\bordermatrix{ 
         & s_{60} & s_{61}  \cr
s_{50} \,& 0,0    & 0,-1    \cr
s_{51} \,& 7,0    & -5,1    \cr},
\end{equation}
(where the $(i,j)$th entry is the pair
$U_5(s_{5i},s_{6j}),U_6(s_{5i},s_{6j})$).
Let $\gamma_1=\gamma_3=1$ and $\gamma_2=\gamma_4=-1$.  Let $U_1$ and $U_2$
be given by
\begin{equation}
\bordermatrix{ 
         & s_{20} & s_{21}  \cr
s_{10} \,& 0,0    & 0,-1    \cr
s_{11} \,& 1,0    & -4,1    \cr},
\end{equation}
and let $U_3$ and $U_4$ be given by
\begin{equation}
\bordermatrix{ 
         & s_{40} & s_{41}  \cr
s_{30} \,& 0,0    & 0,-1    \cr
s_{31} \,& 1,0    & -3,2    \cr}.
\end{equation}
We abbreviate $\sigma_i(s_{i1})$ as $\sigma_i$ by abuse of notation.
At a totally mixed Nash equilibrium $\sigma$ which is hierarchically perfect and
totally mixed with respect to the ENT of Figure \ref{saboteur}, we have
$0=U_5(s_{50},\sigma_6)=U_5(s_{51},\sigma_6)=7(1-\sigma_6)-5\sigma_6=7-12\sigma_6$,
so $\sigma_6=\frac{7}{12}$.  Similarly we have
$0=-(1-\sigma_5)+\sigma_5=2\sigma_5-1$ so $\sigma_5=\frac{1}{2}$.  

We also have
$u_1(s_{10},\sigma_2,\sigma_5,\sigma_6)=U_1(s_{10},\sigma_2)+u_5(\sigma_5,\sigma_6)$,
which we must equate to 
$u_1(s_{11},\sigma_2,\sigma_5,\sigma_6)=U_1(s_{11},\sigma_2)+u_5(\sigma_5,\sigma_6)$,
for hierarchical perfection (here we are ignoring the fact that $\sigma_5$
is a function of $\sigma_1$ and $\sigma_2$).  This gives us that
$0=U_1(s_{10},\sigma_2)=U_1(s_{11},\sigma_2)=(1-\sigma_2)-4\sigma_2=1-5\sigma_2$,
so $\sigma_2=\frac{1}{5}$. Similarly 
we have
$u_2(s_{20},\sigma_1,\sigma_5,\sigma_6)=U_2(s_{20},\sigma_1)-u_5(\sigma_5,\sigma_6)$,
which we must equate to 
$u_2(s_{21},\sigma_1,\sigma_5,\sigma_6)=U_2(s_{21},\sigma_1)-u_5(\sigma_5,\sigma_6)$,
so $U_2(s_{20},\sigma_1)=U_2(s_{21},\sigma_1)$.  This gives
$0=-(1-\sigma_1)+\sigma_1=2\sigma_1-1$, so
$\sigma_1=\frac{1}{2}$. We also have $0=(1-\sigma_4)-3\sigma_4=1-4\sigma_4$, so
$\sigma_4=\frac{1}{4}$, and $0=-(1-\sigma_3)+2\sigma_3=3\sigma_3-1$, so
$\sigma_3=\frac{1}{3}$.  

Finally, we check that
$\sigma_1\sigma_2+(1-\sigma_1)(1-\sigma_2)=\frac{1}{10}+\frac{4}{10}=\frac{1}{2}=\sigma_5$,
and
$\sigma_3\sigma_4+(1-\sigma_3)(1-\sigma_4)=\frac{1}{12}+\frac{6}{12}=\frac{7}{12}=\sigma_6$.
Now given $\sigma_{-1}$, player 1 cannot increase either $U_1$ or $u_5$ by
changing only $\sigma_1$, so player 1 cannot increase $u_1$.  Similarly,
player 2 can neither increase $U_1$ nor {\em decrease} $u_5$ by changing
only $\sigma_2$, so player 2 cannot increase $u_2$.  In this way, we see
that $\sigma$ is a Nash equilibrium of the actual game.

A strategy profile of the actual players is a point in the product of
probability simplices corresponding to their actual strategy spaces.
When we pass to an emergent player one level up, we project the product of
simplices for the actual players below that emergent player to a smaller
dimensional simplex, the space of emergent mixed strategies of this
emergent player.  That we are able to do this means that the payoffs to
other actual players, not below this emergent player, depend only on the
choice of a point in the smaller dimensional simplex by these actual
players.

We can use ENTs to analyze certain cooperative games.  We
consider each coalition to be an emergent player.  An actual player's pure
strategies specify the highest level of coalition to join.  So the number
of its pure strategies is the number of its ancestors in the tree(including
itself).  Each coalition forms if all its descendants agree to join it,
otherwise it doesn't form.  The number of pure strategies of a coalition is
one more than the number of its ancestors (including itself).  Its pure
strategies correspond either to the highest level of coalition containing
this coalition which its members have agreed to form, or to not forming
this coalition itself.  The function $U_v$ for each coalition $v$ is zero
if the coalition forms and is equal to the \emph{value} of the coalition if
it does form; it does not depend on the actions of $v$'s siblings.  The
number $\gamma_{vw}$ represents $v$'s share of the gain from the larger
coalition $w$, if it forms.  

Note that a given ENT does not allow all possible subsets of
players to form coalitions, but only certain ones.  We could extend the
definition to all possible subsets by positing that for any partition of a
coalition into subcoalitions not in the tree, the subcoalitions receives
the same utility by joining or not joining the coalition.  Thus not all
cooperative games correspond to ENTs.  Those that do, however, may often
occur in modeling real situations.
\bibliographystyle{abbrv}
\bibliography{thesis}

\end{document}